\documentclass[12pt,reqno]{amsart}

\usepackage{amsmath}
\usepackage{amssymb}
\usepackage[dvips]{graphicx}

\usepackage{pstricks}
\newtheorem{theo}{Theorem}
 
\newtheorem{prop}[theo]{Proposition}
\newtheorem{coro}[theo]{Corollary}

\theoremstyle{remark}

\begin{document}

\title[Andr\'e permutations, right-to-left and \\ left-to-right minima]%
{Andr\'e permutations, right-to-left and \\ left-to-right minima}
\author{Filippo Disanto*}
\date{}

\thanks{\hbox{\hskip-10pt}*Email: {\tt fdisanto@uni-koeln.de}}

\begin{abstract}
We provide enumerative results concerning right-to-left minima and left-to-right minima in Andr\'e permutations of the first and second kind. For both the two kinds, the distribution of right-to-left and left-to-right minima is the same. We provide generating functions and associated asymptotics results. Our approach is based on the tree-structure of Andr\'e permutations.
\end{abstract}

\maketitle

\thispagestyle{myheadings}
\font\rms=cmr8 
\font\its=cmti8 
\font\bfs=cmbx8

\markright{\its S\'eminaire Lotharingien de
Combinatoire \bfs ?? \rms (201?), Article~????\hfill}
\def\thepage{}

\section{Introduction}

\emph{Andr\'e} permutations have been introduced in \cite{foata} and extensively studied in the literature especially because of their relations with other combinatorial structures~\cite{strehl2, strehl, hetiei, hetiei2, stanleyundici}. For instance, the \emph{cd}-index of the Boolean algebra may be computed by summing the \emph{cd}-variation monomials of Andr\'e permutations~\cite{stanleyundici}. 

It is possible to distinguish among two types of Andr\'e permutations: those of the \emph{first} kind $\mathcal{A}^{(1)}$ and those of the \emph{second} kind $\mathcal{A}^{(2)}$. The two classes are equinumerous. The $n$-th Euler number $e_n = [z^n] \sec(z) + [z^n] \tan(z)$ counts Andr\'e permutations of size $n$. The first terms are $e_0=1,e_1=1,e_2=1,e_3=2,e_4=5,e_5=16,\dots$. Classically, Euler numbers only refer to \emph{secant} numbers, the (even) coefficients of the Taylor expansion of $\sec(z)$. The (odd) coefficients of the Taylor expansion of $\tan(z)$ are called \emph{tangent} numbers. Here, with an abuse of terminology, we take the Euler numbers as the sum of these two coefficients. 
Besides Andr\'e permutations, Euler numbers give the enumeration of several other combinatorial structures. In particular, they also count rooted binary un-ordered \emph{increasing} trees. In \cite{foata} the authors describe two bijections  - denoted here by $\phi_1$ and $\phi_2$ - which maps Andr\'e permutations of both kinds onto this class of trees and viceversa. Based on this correspondence,  two classical permutation statistics, such as right-to-left minima (rlm) and left-to-right minima (lrm), have a natural interpretation in terms of paths of the associated trees. 

In this work we indeed focus on the enumeration of Andr\'e permutations according to the parameters number of right-to-left minima and number of left-to-right minima.  
To the best of our knowledge, these permutation statistics have not been investigated before in this context.

In Section~\ref{inizio}, we show that the statistic number of right-to-left minima has the same distribution on each one of the two sets $\mathcal{A}^{(1)}_n$ and $\mathcal{A}^{(2)}_n$. The same holds for the number of left-to-right minima and, more generally, in the case of their joint distribution. Without loss of generality, we then focus on one type of Andr\'e permutations, those of the second kind $\mathcal{A} = \mathcal{A}^{(2)}$. For the joint enumeration according to right-to-left and left-to-right minima a functional equation for the associated trivariate generating function is provided.

In Section~\ref{L}, we find the bivariate generating function which counts Andr\'e permutations $\mathcal{A}$ with respect to the size and the number of right-to-left minima.  
As a result, fixing the number of right-to-left minima, we provide a combinatorial formula which describes the desired enumeration in terms of Euler numbers. As a corollary to the results of this section we have a correspondance between number of right-to-left minima in Andr\'e permutations and number of cycles in the so-called \emph{cycle-up-down} permutations introduced in \cite{deutsch}. This will need to be further investigated.

In Section~\ref{nelbosco}, we study the number of left-to-right minima. We give a functional equation for the associated bivariate generating function. We show how the number of permutations of size $n+1$ with $2$ left-to-right minima is related to the total number of right-to-left minima in permutations of size $n$. Finally, we study Andr\'e permutations with a generic - but fixed - number of left-to-right minima providing asymptotic estimates. 

\section{Preliminaries} \label{def}

The set of permutations of size $n$ is denoted by $\mathcal{S}_n$. If $\pi = (\pi_1 \pi_2 \dots \pi_n ) \in \mathcal{S}_n$, the set of its \emph{left-to-right minima} is denoted by $\mathrm{lrm}(\pi)$ and its elements are those entries $\pi_i$ such that if $j<i$, then $\pi_i < \pi_j$. We denote by $\mathrm{rlm}(\pi)$ the set of \emph{right-to-left minima} and we remind to the reader that $\pi_i \in \mathrm{rlm}(\pi)$ if $j>i$ implies $\pi_i< \pi_j$.

\pagenumbering{arabic}
\addtocounter{page}{1}
\markboth{\SMALL FILIPPO DISANTO}{\SMALL 
ANDR\'E PERMUTATIONS, RIGHT-TO-LEFT AND LEFT-TO-RIGHT MINIMA}

A binary \emph{increasing} tree is a rooted, \emph{un-ordered} tree with nodes of outdegree $0,1$~or~$2$. Nodes of outdegree $0$ are also called the $leaves$ of the tree. Moreover, for such a tree, we require that each of the $n$ nodes is bijectively labelled by a number in $\{1,2,...,n \}$ in a way that, going from the root to any leaf, we always find an increasing sequence of numbers. If $x$ and $y$ are two nodes, we write $x \prec y$ when the label of $x$ is less than the label of $y$. 
The set of binary increasing trees is denoted by $\mathcal{B}$ while we use the symbol $\mathcal{B}_n$ to denote the subset of $\mathcal{B}$ made of those trees with $n$ nodes.

\begin{figure}
\begin{center}
\includegraphics*[scale=.6,trim=0 0 0 0]{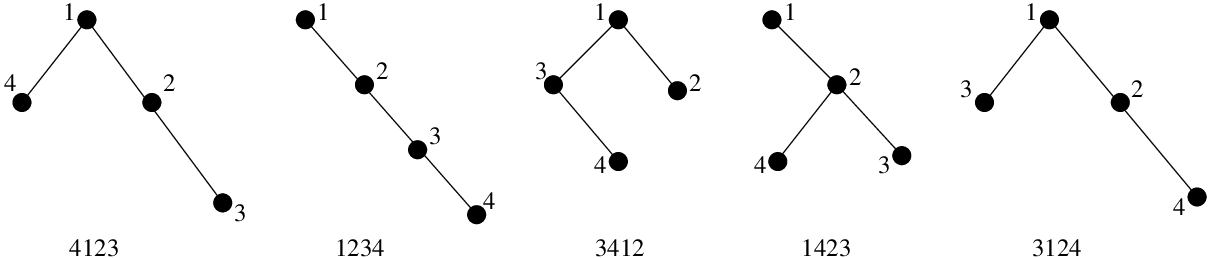}
\end{center}
\caption{The trees in $\mathcal{B}_4$ and the associated Andr\'e permutations of the second kind.}\label{alberetti}
\end{figure}

Observe that each tree in $\mathcal{B}$ can be drawn in the plane in a \emph{unique} way respecting the following two conditions $(A_2)$ and $(B_2)$:

\begin{itemize}
\item[$(A_2)$] if a node has only one child, then this child is drawn on the \emph{right} of its direct ancestor;
\item[$(B_2)$] if a node $x$ has two children $y$ and $z$, let $t_y$ (resp. $t_z$) be the set of nodes in the subtree generated by $y$ (resp. $z$). If $ y=\min(t_y) \prec  z=\min(t_z)$, then $y$ is drawn on the \emph{right} of $x$ while $z$ on the \emph{left}.
\end{itemize} 
In Fig.~\ref{alberetti} we show those trees belonging to $\mathcal{B}_4$ drawn respecting the previous two conditions.

Conditions $(A_2,B_2)$ are not the only possible ones that allow a unique planar representation for each tree in $\mathcal{B}$. Another couple of conditions is for instance

\begin{itemize}
\item[$(A_1)$] if a node has only one child, then this child is drawn on the \emph{right} of its direct ancestor;
\item[$(B_1)$] if a node $x$ has two children $y$ and $z$, let $t_y$ (resp. $t_z$) be the set of nodes in the subtree generated by $y$ (resp. $z$). If $ \max(t_y) \prec  \max(t_z)$, then $z$ is drawn on the \emph{right} of $x$ while $y$ on the \emph{left}.
\end{itemize} 

The sets of \emph{Andr\'e} permutations $\mathcal{A}^{(2)}$ and $\mathcal{A}^{(1)}$ can be defined in several  equivalent ways, see for instance Section~2~of~\cite{hetiei}.  
Since they are both subsets of $\mathcal{S}_n$ equinumerous to $\mathcal{B}_n$, we choose to characterize their permutations according to two injective maps $\phi_2,\phi_1:\mathcal{B}_n \rightarrow \mathcal{S}_n$ (see \cite{foata}). For $\phi_2$ (resp. $\phi_1$) the procedure is:
\begin{itemize}
\item[(1)] given $t \in \mathcal{B}_n$, draw $t$ according to $(A_2,B_2)$ (resp. $(A_1,B_1)$);
\item[(2)] each leaf collapses into its direct ancestor whose label is then modified receiving on the left the label of the left child (if any) and on the right the label of its right child. We obtain in this way a new tree whose nodes are labelled with sequences of numbers; 
\item[(3)] starting from the obtained tree go to step (2).
\end{itemize}

The algorithms $\phi_2, \phi_1$ end when the tree $t$ is reduced to a single node whose label is then a permutation $\phi_2(t), \phi_1(t)$ of size $n$. Note that, without considering step (1) but only (2) and (3), the procedures give a well-known \cite{stan} bijection $\psi$ between \emph{ordered} binary increasing trees $\tilde{\mathcal{B}}_n$ and the entire set of permutations of size $n$. 

\smallskip

The sets $\mathcal{A}_n^{(i)}$  can be defined as $\mathcal{A}_n^{(i)} = \{ \phi_i(t) \in \mathcal{S}_n : t \in \mathcal{B}_n \}$ (with $i = 2,1$). Looking at Fig.~\ref{alberetti}, the corresponding permutations in $\mathcal{A}_4^{(2)}$ are (from left to right) $(4123), (1234), (3412), (1423)$ and $(3124)$. For the same size $n=4$, the permutations in $\mathcal{A}_4^{(1)}$ are $(2314), (1234), (2134), (1324)$ and $(3124)$.

\smallskip

An equivalent definition of Andr\'e permutations can be given in terms of the so-called $x$\emph{-factorizations} of permutations, see Definition~1 and Definition~2 of \cite{hetiei2} and the related references. The equivalence is easily recovered by observing that - following notations of \cite{hetiei2} - the $\lambda$-part of the $x$-factorization of a permutation $\pi$ corresponds to the left-subtree of the node $x$ in the \emph{ordered} binary increasing tree $\psi^{-1}(\pi)$. Similarly, the $\rho$-part of the $x$-factorization corresponds to the right sub-tree of $x$ in $\psi^{-1}(\pi)$.     

Andr\'e permutations, as binary increasing trees, are enumerated, with respect to the size, by the so called \emph{Euler} numbers $(e_n)_{n\geq 0}$ whose exponential generating function satisfies

$$\int E^2 = 2E-z-2 \footnote{\text{We will often adopt the notation} $\int f(z) = \int_0^{z} f(a) da$.}$$ and therefore is equal to

$$E(z)= \sec(z) + \tan(z).$$

The first terms of the sequence are: $1,1,1,2,5,16,61,272,1385,...$ and they correspond to entry $A000111$ in \cite{sloane}. Furthermore, expanding $E(z)$ near the dominant singularity $z=\pi/2$, we easily recover an asymptotic approximation for the coefficients 
\begin{equation}\label{brioscia}
\frac{e_n}{n!} \sim \frac{4}{\pi} \left( \frac{2}{\pi}  \right)^n.
\end{equation}

\section{Enumeration of right-to-left minima and left-to-right minima}

In this section we study enumerative properties of right-to-left minima and left-to-right minima in Andr\'e permutations. In Section~\ref{inizio}, these statistics are jointly studied. In Section~\ref{L}, we focus on the number of right-to-left minima, while, in Section~\ref{nelbosco}, we investigate left-to-right minima.

\subsection{Joint enumeration}\label{inizio}
 
Through the bijection $\psi: \tilde{\mathcal{B}}_n \rightarrow \mathcal{S}_n$ described in Section~\ref{def}, one can see that, for any given pemutation $\pi$, the set $\mathrm{rlm}(\pi)$ corresponds to the nodes visited in the tree $\psi^{-1}(\pi)$ starting from the root and performing only right-steps. Similarly, the set $\mathrm{lrm}(\pi)$ corresponds to the nodes visited in the tree $\psi^{-1}(\pi)$ starting from the root and performing only left-steps.

\begin{figure}
\begin{center}
\includegraphics*[scale=1.0,trim=0 0 0 0]{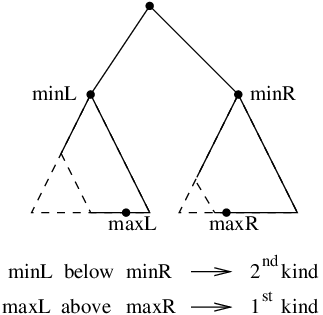}
\end{center}
\caption{Recursive decomposition of $t_i = \phi_i^{-1}(\pi_i)$, with $\pi_i \in \mathcal{A}^{(i)}$ ($i = 2,1$). The $\min$-node and $\max$-node of each root sub-tree are highlighted. The left corner of these subtrees can be empty according to the right orientation of single nodes (conditions $A_i$).}\label{alba}
\end{figure}

Let $\pi_2 \in \mathcal{A}_n^{(2)}$ and $\pi_1 \in \mathcal{A}_n^{(1)}$, consider $t_2 = \phi_2^{-1}(\pi_2)$ and $t_1 = \phi_1^{-1}(\pi_1)$. If $n>1$ then, for $i=2,1$, the tree $t_i$ consists of two trees, $t_{i,\text{left}}$ and $t_{i,\text{right}}$, appended to its root on the left and on the right respectively. Clearly, $\mathrm{rlm}(\phi_i^{-1}(t_i)) = 1 +  \mathrm{rlm}(\phi_i^{-1}(t_{i,\text{right}})) $ and $\mathrm{lrm}(\phi_i^{-1}(t_i)) = 1 +  \mathrm{lrm}(\phi_i^{-1}(t_{i,\text{left}})) $, see Fig.~\ref{alba}. 

Furtermore, observe that, in both cases $i=2,1$, there are exactly $${{|t_{i,\text{left}}| + |t_{i,\text{right}}| - 1}\choose{|t_{i,\text{left}}|}}$$ ways of merging the ranking of $t_{i,\text{left}}$ with the ranking of $t_{i,\text{right}}$ that create a tree drawn according to conditions ($A_i, B_i$). When $i=2$, we have to put the root of $t_{i,\text{right}}$ \emph{above} the root of $t_{i,\text{left}}$  while, when $i=1$, we put the $\max$-node of $t_{i,\text{right}}$ \emph{below} the $\max$-node of $t_{i,\text{left}}$ (the $\max$-node is always a leaf). Also note that, when $|t_{i,\text{left}}|=0$, the previous binomial expression returns $1$.

From this considerations, it follows that, from an enumerative point of view, the same recursive construction describes the distribution of right-to-left minima and left-to-right minima in Andr\'e permutations of the first and second kind.  

Without loss of generality, we decide to focus on Andr\'e permutations of the second kind. We thus set $\mathcal{A} = \mathcal{A}^{(2)}$, $\phi= \phi_2$ and, if not specified otherwise, we draw each tree $t \in \mathcal{B}$ according to $(A_2,B_2)$.

\smallskip

The exponential generating function $$H = H(x,y,z) = \sum_{\pi \in \mathcal{A}} \frac{x^r y^l z^n}{n!},$$ where $r=|\mathrm{rlm}(\pi)|, l=|\mathrm{lrm}(\pi)|$ and $n=\text{size}(\pi)$, satisfies the functional equation

$$H = 1 + xyz + \sum_{\pi_1=t_{\text{right}} \neq \emptyset} \sum_{\pi_2=t_{\text{left}}} x^{r_1 + 1} y^{l_2 + 1}\frac{z^{n_1+n_2+1}}{(n_1+n_2+1)!} \cdot {{n_1+n_2-1}\choose{n_2}}.$$

Taking twice the derivative respect to $z$ we obtain
\begin{equation}\label{nonni}
\frac{\partial^2 H }{\partial z^2} = x y \frac{\partial H(x,1,z)}{\partial z} \, H(1,y,z)
\end{equation}
which gives 
\begin{equation}\label{tagliaerba}
H = 1 + xyz + xy \int \int \frac{\partial H(x,1,z)}{\partial z} \, H(1,y,z) dz dz.
\end{equation}

Equation (\ref{tagliaerba}) can be used recursively to compute the polynomials $H_i(x,y)= \sum_{\pi \in \mathcal{A}_i} x^r y^l$. When $0\leq i \leq 5$ we have
\begin{eqnarray}\nonumber
H_0 &=& 1 ;\\\nonumber 
H_1 &=& xy ; \\\nonumber 
H_2 &=& x^2 y ; \\\nonumber
H_3 &=& x^3 y + x^2 y^2 ; \\\nonumber
H_4 &=& x^4 y + 2 x^3 y^2 + x^3 y + x^2 y^2 ; \\\nonumber
H_5 &=&  x^5 y + 3 x^4 y^2 + 3 x^4 y + 6 x^3 y^2 +  x^3 y +  x^2 y^3 + x^2 y^2. \nonumber  
\end{eqnarray} 

Furthermore, considering that $H(1,1,z) = E(z)$ and that $E'(z) = \frac{1}{1-\sin(z)}$, equation (\ref{nonni}) becomes
\begin{equation}\label{rait}
\frac{\partial^2 H(x,1,z) }{\partial z^2} = x \frac{\partial H(x,1,z)}{\partial z} \, E(z)
\end{equation}

when we consider $y=1$, while it gives 

\begin{equation}\label{lefttu}
\frac{\partial^2 H(1,y,z) }{\partial z^2} = y E'(z) \, H(1,y,z)
\end{equation}

when we take $x=1$. In the following sections we will study (\ref{rait}) and (\ref{lefttu}) as they respectively provide the enumeration of Andr\'e permutations with respect to the number of right-to-left minima and left-to-right minima.

\subsection{Right-to-left minima} \label{L}

Here we focus on \emph{right-to-left minima} statistic using the symbol $\mathcal{A}^{R}_{n,r}$ to denote the subset of $\mathcal{A}_n$ made of those permutations $\pi$ with $|\mathrm{rlm}(\pi)|=r$.

Defining 
\begin{equation} \nonumber
F(x,z) = \left( \frac{1}{1-\sin(z)} \right)^x, 
\end{equation}
it is easy to check that 
$$\frac{\partial F(x,z)}{\partial z} = x F(x,z) \cdot E(z).$$

Thus, setting $$F(x,z) = \frac{1}{x} \frac{\partial H(x,1,z)}{\partial z},$$ we have that $H(x,1,z)$ satisfies (\ref{rait}). Considering $F(1,z)$ provides the (shifted) exponential generating function for Euler numbers.
In other words, we have the following result

\begin{prop}\label{collara}
The (shifted) exponential generating function counting Andr\'e permutations with respect to the size $n$ and number of righ-to-left minima $r$ is given by
\begin{equation} \nonumber
F(x,z) = \left( \frac{1}{1-\sin(z)} \right)^x = \sum_{\pi \in \mathcal{A}} \frac{ x^{r-1} z^{n-1} }{(n-1) \, !} . \nonumber
\end{equation}
\end{prop}

\bigskip

The first terms of $|\mathcal{A}^{R}_{n,r}|$ are thus given by the following table.

\begin{center} 
\begin{tabular}{|c|ccccccccc|}
\hline
n/r & 2 & 3 & 4 & 5 & 6 & 7 & 8 & 9 & 10 \\\hline
2 & 1 & 0 & 0 & 0 & 0 & 0 & 0 & 0 & 0  \\ 
3 & 1 & 1 & 0 & 0 & 0 & 0 & 0 & 0 & 0 \\ 
4 & 1 & 3 & 1 & 0 & 0 & 0 & 0 & 0 & 0 \\ 
5 & 2 & 7 & 6 & 1 & 0 & 0 & 0 & 0 & 0 \\ 
6 & 5 & 20 & 25 & 10 & 1 & 0 & 0 & 0 & 0 \\ 
7 & 16 & 70 & 105 & 65 & 15 & 1 & 0 & 0 & 0\\ 
8 & 61 & 287 & 490 & 385 & 140 & 21 & 1 & 0 & 0\\ 
9 & 272 & 1356 & 2548 & 2345 & 1120 & 266 & 28 & 1 & 0   \\  
10 & 1385 & 7248 & 14698 & 15204 & 8715 & 2772 & 462 & 36 & 1 \\\hline
\end{tabular} 
\end{center}

Note that Euler numbers are the entries of the first column. 
Furthermore observe that, looking at the table column by column, one has

$$\left( \frac{\partial ^{r} \, F}{\partial x ^{r}} \right)_{x=0}=\bigg[-\ln \big( 1-\sin(z) \big)\bigg]^r$$ 

and then

\begin{equation}\label{rota}
\frac{1}{r!}\bigg[-\ln \big( 1-\sin(z) \big)\bigg]^r = \sum_{\pi, |\mathrm{rlm}(\pi)|=r+1} \frac{z^{n-1} }{(n-1) \, !}. 
\end{equation}

Given that $\int E(z) = -\ln \big( 1-\sin(z) \big),$ as a corollary we have 

\begin{prop}\label{ele}
For every fixed $r \geq 1$
\begin{equation}\label{balotelli}
\frac{1}{r!}\left[ \sum_{n\geq 1} \frac{e_{n-1}}{n!}z^n \right]^{r}=\sum_{n\geq l}\frac{|\mathcal{A}^{R}_{n+1,r+1}|}{n!} z^n,
\end{equation}
where  
$e_0=1, e_1=1, e_2=1, e_3=2, e_4=5, e_5=16, e_6=61,...$
are Euler numbers.
\end{prop}

For a fixed $r\geq 2$, the asymptotic behaviour of the sequence $|\mathcal{A}^{R}_{n,r}|$ can also be examined at this point. Observe that near the dominant singularity $z=\pi/2$, from 
$$-\ln\big(1-\sin(z)\big) = \ln(2) - 2 \ln(z-\pi/2) + 1/12 (z-\pi/2)^2 + \mathcal{O}\left( (z-\pi/2)^4  \right),$$ 
we have the following approximation
\begin{eqnarray}\nonumber
\bigg[-\ln\big(1-\sin(z)\big) \bigg]^r &=& \bigg[ \ln(2) - 2 \ln(z-\pi/2)  \bigg]^r + \mathcal{O}\left( z-\pi/2 \right) \\\label{paguro} 
&=& 2^r \bigg[-\ln(z-\pi/2)  \bigg]^r + 2^{r-1} r \ln(2)\bigg[-\ln(z-\pi/2)  \bigg]^{r-1} \\\nonumber
&& + \mathcal{O}\left( \bigg[-\ln(z-\pi/2)  \bigg]^{r-2}  \right). 
\end{eqnarray}

Rewriting $\ln(z- \pi/2) = \ln(-\pi/2) + \ln\left(1-\frac{2z}{\pi} \right)$, by Th.VI.2 of \cite{ancomb} (see special cases formula (27)) we have 
$$[z^n]\bigg[-\ln(z-\pi/2)  \bigg]^r \sim ( 2/\pi)^n n^{-1}\left[ C_1 \big(\ln(n)\big)^{r-1} + \mathcal{O}\left( \big(\ln(n)\big)^{r-2}  \right)  \right]$$
and similarly
$$[z^n]\bigg[-\ln(z-\pi/2)  \bigg]^{r-1} \sim ( 2/\pi)^n n^{-1}\left[ C_2 \big(\ln(n)\big)^{r-2} + \mathcal{O}\left( \big(\ln(n)\big)^{r-3}  \right)  \right],$$
where $C_1,C_2$ are positive constants.

Furthermore, by using Th.VI.3 \cite{ancomb} for the $\mathcal{O}$-transfer, we have that 
$$[z^n]\bigg[ \mathcal{O}\left( \bigg[-\ln(z-\pi/2)  \bigg]^{r-2}  \right)  \bigg] = \mathcal{O}\left( (2/\pi)^n n^{-1}  \bigg( \ln(n) \bigg)^{r-2} \right).$$

Finally, by applying Th.VI.4 \cite{ancomb} to (\ref{paguro}) and recalling (\ref{rota}),
we obtain the following result.
\begin{prop}
For a fixed $r \geq 1$ and $n \rightarrow \infty$, we have the asymptotic equivalence 
\begin{equation}
\frac{|\mathcal{A}^{R}_{n+1,r+1}|}{n!} = [z^n]\bigg[-\ln\big(1-\sin(z)\big) \bigg]^r \sim k_r \cdot n^{-1} \left(\frac{2}{\pi}\right)^n \bigg( \ln(n) \bigg)^{r-1},
\end{equation}
where $k_r$ is a positive constant depending on $r$.
\end{prop}

\bigskip

We conclude this section recalling that in Chapter 7 of \cite{johnson} the author studies a family of polynomials corresponding to the rows of the previous table. He also shows a criterion according to which each row defines a partition of the set of \emph{up-down} permutations of a given size. Furthermore, in \cite{deutsch} the authors prove that the rows of the previous table also provide the enumeration of the so called \emph{cycle-up-down} permutations with respect to the size and to the number of cycles. It is then natural to ask for a bijection between the permutations in $\mathcal{A}_{n+1}$ and the cycle-up-down ones of size $n$ enlightening the correspondence between right-to-left minima and cycles.

\subsection{Left-to-right minima}\label{nelbosco}

In the previous section we have enumerated the permutations in $\mathcal{A}$ with respect to the size and to the number of right-to-left minima. 
Here we study the cardinality of $\mathcal{A}^{L}_{n,l}$, that is the subset of $\mathcal{A}_n$ with $|\mathrm{lrm}(\pi)|=l$. 

Using the polynomials $H_i$ of Section~\ref{inizio}, we have computed the entries of the following table showing, for all $(n,l)\in \{1,...,10 \}^2$, the number of permutations in $\mathcal{A}^{L}_{n,l}$.

\begin{center} 
\begin{tabular}{|c|cccccccccc|}
\hline
n/l & 1 & 2 & 3 & 4 & 5 & 6 & 7 & 8 & 9 & 10 \\\hline
1 &  1 & 0  & 0  & 0  & 0  & 0  & 0  & 0  & 0  & 0 \\
2 &  1 & 0  & 0  & 0  & 0  & 0  & 0  & 0  & 0  & 0 \\ 
3 &  1 & 1  & 0  & 0  & 0  & 0  & 0  & 0  & 0  & 0 \\ 
4 &  2 & 3  & 0  & 0  & 0  & 0  & 0  & 0  & 0  & 0 \\ 
5 &  5 & 10  & 1   & 0  & 0  & 0  & 0  & 0  & 0  & 0 \\ 
6 &  16 & 38  & 7  & 0  & 0  & 0  & 0  & 0  & 0  & 0 \\ 
7 & 61  & 165  & 45  & 1  & 0  & 0  & 0  & 0  & 0  & 0 \\ 
8 & 272  & 812  & 288  & 13  & 0  & 0  & 0  & 0  &  0 & 0 \\ 
9 &  1385 & 4478  & 1936  & 136  & 1  & 0  & 0  & 0  & 0  &  0  \\  
10 & 7936 &  27408 & 13836  & 1320  & 21  & 0  &  0 & 0  & 0  & 0  \\\hline
\end{tabular} 
\end{center}

In the first column we find (shifted) Euler numbers. It is also interesting to observe that the entries in the second column
$$
1,3,10,38,165,812,4478,27408,184529,1356256,10809786,92892928...
$$

belong to sequence A186367 of \cite{sloane}. This sequence counts the number of cycles in all \emph{cycle-up-down} permutations of size $n$ (see also \cite{deutsch}) and, furthermore, it is strongly related to the total number of right-to-left minima in the permutations of $\mathcal{A}$ having fixed size. Indeed, we will prove that the exponential generating function associated with the non-zero entries of the column $r=2$ of the table above is given by 

$$\left( \frac{\partial F}{\partial x} \right)_{x=1} =\frac{-\ln \big(1-\sin(z)\big)}{1-\sin(z)},$$
where $F$ is the same of Proposition~\ref{collara}. 

\begin{figure}
\begin{center}
\includegraphics*[scale=.5,trim=0 0 0 0]{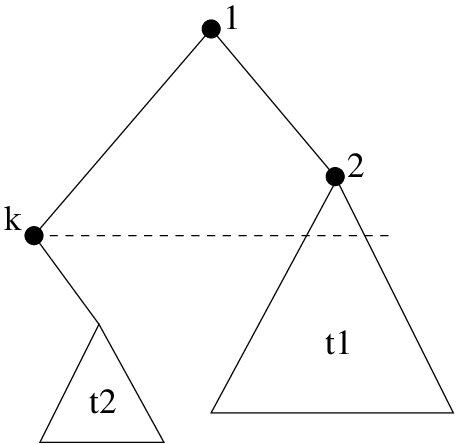}
\end{center}
\caption{Decomposition of a tree with two left-to-right minima.}\label{jollo}
\end{figure}

In order to prove the correspondence, we observe that each tree $\phi^{-1}(\pi)$, such that $|\mathrm{lrm}(\pi)=2|$, can be decomposed as shown in Fig.~\ref{jollo}. 

In particular, note that tree $t_1$ must contain at least one node (labelled with $2$) while tree $t_2$ could be empty. The class of trees consisting of a (possibly empty) tree appended to a node - denoted by $k$ in Fig.~\ref{jollo} - is counted by the exponential generating function 
\begin{eqnarray}\nonumber
f_{t_2}(z) &=& \sum_{m>0} \frac{\tilde{e}_m z^m}{m!}=\int E(z)=-\ln(1-\sin(z)) \,\,\, (\mathrm{with}\,\, \tilde{e}_m=e_{m-1}) \\\nonumber
\end{eqnarray}

while 

\begin{eqnarray} \nonumber
f_{t_1}(z) &=& \sum_{n>0} \frac{e_n z^n}{n!}=E(z)-1 \\\nonumber 
\end{eqnarray}
counts those trees having at least one node. Appending $t_1$ of size $n$ and $t_2$ of size $m-1$ as shown in Fig.~\ref{jollo}, we can build exactly ${{n+m-1}\choose{m}}$ different trees. 
It follows that, in the previous table, the entries $n\geq 1$ of the column $r=2$  correspond to the coefficients of the following exponential generating function 
\begin{equation}\nonumber
g_2(z)=\sum_{n>0} \sum_{m>0} \frac{e_n \tilde{e}_m z^{n+m+1}}{(n+m)(n+m+1)(m!)(n-1)!}.
\end{equation}

Finally observe that

\begin{equation}\nonumber
g_2'' = f_{t_2} \cdot f_{t_1}' = \ln \left( \frac{1}{1-\sin(z)} \right) \cdot \left(\frac{1}{1-\sin(z)} \right)=\left( \frac{\partial F}{\partial x} \right)_{x=1}.
\end{equation}

Given the above calculations, we obtain the following result

\begin{prop}
The following equality holds for all $n \geq 2$:
\begin{eqnarray}\nonumber
|\mathcal{A}^{L}_{n+1,2}|&=& \sum_{r\geq 2} (r-1)\cdot |\mathcal{A}^{R}_{n,r}|\\\nonumber
&=&(n-1)!\cdot [z^{n-1}]\left(\frac{-\ln \big(1-\sin(z)\big)}{1-\sin(z)}\right) \\\nonumber
\end{eqnarray} 
\end{prop}

from which we have the next corollary 

\begin{coro}\label{rosario}
For $n\geq 2$: 
\begin{equation}\label{media}
|\mathcal{A}^{L}_{n+1,2}| +|\mathcal{A}_n| = \sum_{r\geq 2} r \cdot |\mathcal{A}^{R}_{n,r}|
\end{equation} 
and therefore the expected number of right-to-left minima in a random permutation of $\mathcal{A}_n$ is given by  $1 + | \mathcal{A}^{L}_{n+1,2}|/|\mathcal{A}_n|.$
\end{coro}

\subsubsection{Fixing the number of left-to-right minima}

It is interesting to investigate more in details what happens when we fix the number of left-to-right minima in $\mathcal{A}_n$.
Let $$G_{l}(z)=\sum_{n} \frac{|\mathcal{A}^{L}_{n,l}|}{n!}\cdot z^n$$ and
$$G(y,z) = \sum_{l \geq 0} y^l G_l(z).$$ 

Thus $G = H(1,y,z)$ and, from (\ref{lefttu}), we can write that
\begin{equation}\label{checifaccio}
\frac{\partial ^2 G}{\partial z^2} = y E'(z) \cdot G,
\end{equation}
where $E'(z)=\left( \frac{1}{1-\sin(z)}\right)$.

From (\ref{checifaccio}) we have a recursion for $G_l$.

\begin{prop} \label{ecci}
The family of generating functions $(G_l)_l$
satisfies
\begin{equation}\label{recg}
G_l(z) = \int \int G_{l-1}(z) \cdot E'(z)  
\end{equation}
being $E'(z)=\left( \frac{1}{1-\sin(z)}  \right)$ and $G_1(z)=\int E(z)  = -\ln\big(1-\sin(z)\big)$. 
\end{prop}

Unfortunately equation~(\ref{checifaccio}) does not give an explicit solution for $G$. Still, as we will see later, it can be used to explore the structure of the solution in a neighbourhood of the singularity $z=\pi/2$.

Let us now focus on the exact computation of $G_l$. To do so, one can apply the result of Proposition~\ref{ecci} together with the fact that $E=E(z)$ satisfies $\int E^2= 2E-z-2$. Here we compute explicitly the generating functions $G_l$ for the first values of $l$, say $l=1,2,3$, enlightening the correspondance with the generating function for Euler numbers. If we define $$\int^{(i)} f = \stackrel{i \mathrm{-times}}{\overbrace{\int\int\dots\int} f},$$ 
for $l=1,2$ we have
\begin{align}\nonumber
G_1 =& \int E \\\nonumber
G_2 =& \left(\int^{(2)} \left(\int E\right) E'\right) = \left(\int \left(E \int E \right) \right) -\int^{(2)} E^2 \\\nonumber 
=& \frac{1}{2}\cdot \left( \int E \right)^2 - \int^{}\left( 2E-z-2\right) \\\nonumber
=& \frac{1}{2}\cdot \left( \int E \right)^2 -2 \left(\int E\right) + \frac{z^2}{2} + 2z, \\\nonumber
\end{align}

while, for $l=3$, we obtain:

\begin{align}\nonumber
G_3 =& \left(\int^{(2)} \frac{(\int E)^2}{2} E'\right) -2 \left(\int^{(2)} \left( \int E \right) E' \right) + \left(\int^{(2)} \left(\frac{z^2}{2} + 2z\right)E' \right)\\\nonumber
=& \left( \int E \frac{(\int E)^2}{2} \right) - \left( \int^{(2)} E^2 \left(\int E\right) \right)-2\left[ \left(\int E \left( \int E \right) \right) -\int^{(2)} E^2  \right]\\\nonumber
& + \left(\frac{z^2}{2}+2z \right)\left(\int E \right)+ \left(-2z-4 \right)\left(\int^{(2)} E \right)+ 3\left(\int^{(3)} E \right)\\\nonumber
=& \frac{1}{6}\cdot \left( \int E \right)^3 -\left[ \left(\int (2E-z-2)\left(\int E \right) \right)-\left(  \int^{(2)} (2E-z-2)E\right)  \right]\\\nonumber
& -2\left[ \frac{1}{2}\cdot \left(\int E \right)^2 -\int (2E-z-2)\right]+ \left(\frac{z^2}{2}+2z \right)\left(\int E \right)\\\nonumber
& + \left(-2z-4 \right)\left(\int^{(2)} E \right)+ 3\left(\int^{(3)} E \right)\\\nonumber
=& \frac{1}{6} \left( \int E\right)^3 -2\left( \int E\right)^2 +\left( 8+2z +\frac{z^2}{2}\right)\left( \int E\right)-2z^2-8z\\\nonumber
& + (-2z -4) \left( \int^{(2)} E \right) +4\left( \int^{(3)}E \right).
\end{align}

Recalling that $$[z^n]\left(\int^{(i)} E(z) \right) = \frac{e_{n-i}}{n!},$$ the previous calculations express $|\mathcal{A}^{L}_{n,l}|$ ($l=1,2,3$) in terms of Euler numbers $e_n$.

For values of $l$ greater than $3$ the computation of $G_l$ becomes more difficult. In this cases, we can still use the results of Proposition~\ref{ecci} to obtain asymptotic estimates of the coefficients $[z^n]G_l(z)$. Using standard methods of analytic combinatorics (see \cite{ancomb}), it is sufficient to know an approximation of the function $G_l$ near its dominant singularity to describe the behaviour of $[z^n]G_l(z)$ for $n \rightarrow \infty$. In this case, the idea is to iteratively recover an approximation for $G_{l+1}$ by integration of an approximation for $(G_l \cdot E')$. 

Near the dominant singularity $z=\pi/2$ we have 
\begin{equation}\label{seno}
E'(z)=\frac{1}{1-\sin(z)}= \frac{2}{\left(\frac{\pi}{2}-z \right)^2} + \mathcal{O}\left(1\right)
\end{equation}
and, for every $A>0$,
\begin{equation}\label{castrellina}
G_1= \int E(z) =\ln\left(\frac{1}{1-\sin(z)} \right)= -2\ln\left(\frac{\pi}{2}-z\right)+\mathcal{O}\left(1\right)
= \mathcal{O}\left( \left( \frac{\pi}{2} - z \right)^{-A} \right).
\end{equation}

Then, as a first approximation, one has
$$(G_1 \cdot E')(z) = \mathcal{O}\left(\left(\frac{\pi}{2}-z\right)^{-2-A}\right),$$ which gives by Proposition~\ref{ecci} and Th.~VI.9 \cite{ancomb} (see case (i))
$$G_2(z) = \mathcal{O}\left(\left(\frac{\pi}{2}-z \right)^{-A} \right).$$
We remark that, by the mentioned theorem, we can obtain a singular approximation of $G_2$ by integrating, according to classical rules, the singular expansion of $(G_1\cdot E')$.  
Iterating the procedure one has that, independently on $l$, for every $A>0$ 
\begin{equation}\label{titti}
G_l(z) =  \mathcal{O}\left(\left(\frac{\pi}{2}-z \right)^{-A} \right).
\end{equation}

Applying Th.~VI.3 of \cite{ancomb} to (\ref{titti}) gives the following bound.

\begin{prop}
 When $n$ is large, for every $A>0$ and independently on $l$, we have
\begin{equation}\label{proba}
\frac{|\mathcal{A}^{L}_{n,l}|}{n!}=[z^n]G_l(z)=  \, \mathcal{O}\left( \left(\frac{2}{\pi}\right)^n\cdot n^{A-1} \right).
\end{equation} 
\end{prop}

Recalling that $\frac{|\mathcal{A}_{n}|}{n!} \sim \frac{4}{\pi} \left( \frac{2}{\pi} \right)^n$ (see (\ref{brioscia})), equation (\ref{proba}) gives a measure of how strong is the effect of fixing the number of left-to-right minima in Andr\'e permutations.   

\smallskip

\textbf{Structural properties of $G$ near the singularity.} To conclude our asymptotic analysis we go back to equation~(\ref{checifaccio}) to describe a structural property of the solution $G$. Indeed, treating $y$ as a constant, we can apply Th.~VII.9 of \cite{ancomb} finding that near the \emph{regular} singular point $z=\pi/2$ the desired solution $G$ can be expressed as 
\begin{equation}\nonumber
G = a_{y}\cdot  \left(\frac{\pi}{2}-z \right)^{\frac{1+\sqrt{1+8y}}{2}} A_{y}\left(z-\frac{\pi}{2}\right) +{b}_{y}\cdot \left(\frac{\pi}{2}-z \right)^{\frac{1-\sqrt{1+8y}}{2}} B_{y}\left(z-\frac{\pi}{2}\right),
\end{equation}
where $y$ could in principle appear in $a_{y},A_{y}(z),b_{y},B_{y}(z)$ and the functions $A_{y}(z), B_{y}(z)$ are analytic at $z=0$.

It is interesting to note that, taking $a_{y}=0$ and $b_{y}=B_{y}=1$, one obtains
$$G_{\alpha}=\left(\frac{\pi}{2}-z \right)^{\frac{1-\sqrt{1+8y}}{2}}$$
whose expansion at $y=0$ looks as
\begin{align}\nonumber
G_{\alpha}=& 1-2y \ln(\pi/2 - z) + y^2 \left(4 \ln(\pi/2 - z) + 2 \left[\ln(\pi/2 - z)\right]^2 \right) \\\nonumber
& + y^3\left(-16\ln(\pi/2 - z) -8\left[\ln(\pi/2 - z)\right]^2 -\frac{4}{3}\left[\ln(\pi/2 - z) \right]^3 \right) + \cdots .\nonumber
\end{align}
Based on the approximation for $\int E$ given in (\ref{castrellina}), this reflects the asymptotic behaviour of the expressions for $G_1, G_2$ and $G_3$ which have been previously computed. 
This can be justified observing that $G_{\alpha}$ satisfies 
\begin{equation}\nonumber
\frac{\partial ^2 G_{\alpha}}{\partial z^2} = y \cdot \frac{2}{(\pi/2 -z)^2} \cdot G_{\alpha}
\end{equation}
which is obtained by substituting in (\ref{checifaccio}), i.e. the defining equation for $G$, the term $E'(z)$ by $2/(\pi/2-z)^2$, the latter being the main part of the singular approximation (\ref{seno}). 

\section*{Acknowledgement}

This work was financially supported by grant DFG-SPP1590 from the German Research Foundation.


\begin{thebibliography}{12}


\bibitem{deutsch}
E. Deutsch, S. Elizalde, \emph{Cycle-up-down permutations}, Australas. J. Combin., 50 (2011) 187-199

\bibitem{fill}
J. A. Fill, P. Flajolet, N. Kapur, \emph{Singularity Analysis, Hadamard Products and Tree Recurrences}, J. Comput. Appl. Math. 174 (2005) 271-313.

\bibitem{ancomb}
P. Flajolet, R. Sedgewick, \emph{Analytic Combinatorics}, Cambridge University press, (2009).

\bibitem{foata}
D. Foata, M. P. Sch\"utzenberger, \emph{Nombres d'Euler et permutations alternantes}, unabriged version, 71 pages, University of Florida, Gainesville (1971), available at \texttt{http://www.mat.univie.ac.at/$\sim$slc/}.


\bibitem{strehl2}
D. Foata, V. Strehl, \emph{Rearrangements of the symmetric group and enumerative properties of the tangent and secant numbers}, Math. Z., 137 (1974) 257-264. 

\bibitem{strehl}
D. Foata, V. Strehl, \emph{Euler numbers and variations of permutations}. In: Colloquio Internazionale sulle Teorie Combinatorie, 1973, Tome I (Atti Dei Convegni Lincei 17, 119 131), Accademia Nazionale dei Lincei, Rome (1976).

\bibitem{hetiei}
G. Hetyei, \emph{On the cd-Variation Polynomials of Andr\'e and Simsun Permutations}, Discrete Comput. Geom., 16 (1996) 259-275.

\bibitem{hetiei2}
G. Hetyei, E. Reiner, \emph{Permutation Trees and Variation Statistics}, Europ. J. Combinatorics, 19 (1998) 847-866.


\bibitem{johnson}
W. P. Johnson, \emph{Some polynomials associated with up-down permutations}, Discrete Mathematics, 210 (2000) 117-136.

\bibitem{sloane}
N. J. A. Sloane, \emph{The On-Line Encyclopedia of Integer Sequences}, available at: http://oeis.org/.

\bibitem{stan}
R.P. Stanley, \emph{Enumerative Combinatorics}, vol.1, Wadsworth \& Brooks, Monterey, CA (1986).

\bibitem{stanleyundici}
R. P. Stanley, \emph{Flag f-vectors and the cd-index}, Math. Z., 216 (1994) 483-499.


\end{thebibliography}
\end{document}